\documentclass{article}

\usepackage[tbtags]{amsmath}
\usepackage{amsfonts,amssymb,amscd,graphicx}
\usepackage{mathrsfs}
\usepackage{hypbmsec}

\numberwithin{equation}{section}

\newcommand\qed{\hfill$\Box$}
\newcommand\op\operatorname
\def\bcdot{\,\boldsymbol\cdot\,}

\newenvironment{pf}[1]{\ifdim\lastskip<2.1pt\vskip%
-\lastskip\else\fi\vskip2pt\textbf{#1.}}{\vskip2pt}

\def\itemi{\par\parindent13.8pt{\textup{(i)}}\kern4pt}
\def\itemii{\par\parindent11pt{\textup{(ii)}}\kern4pt}
\def\itemiii{\par\parindent8pt{\textup{(iii)}}\kern4pt}
\def\itemiv{\par\parindent8.5pt{\textup{(iv)}}\kern4pt}

\renewcommand\epsilon{\varepsilon}

\newtheorem{proposition}{Proposition}
\newtheorem{theorem}{Theorem}
\newtheorem{lemma}{Lemma}
\newtheorem{corollary}{Corollary}

\newtheorem{definition}{Definition}

\begin{document}

\title{Combinatorial encoding of Bernoulli schemes and the asymptotic behavior of Young tableaux}

\author{A.\,\,M.~Vershik\thanks{%
St.~Petersburg Department of Steklov Institute of Mathematics and
St.~Petersburg State University, St.~Petersburg, Russia; Institute for Information Transmission Problems, Moscow, Russia.
E-mail: {\tt avershik@gmail.com}. Supported by the RSF grant 17-71-20153.}}

\date{}

\maketitle

\begin{abstract}
We consider two examples of a fully decodable combinatorial encoding of Bernoulli schemes: the encoding via Weyl simplices and the  much more complicated encoding via the  RSK (Robinson--Schensted--Knuth) correspondence. In the first case, the decodability  is a quite simple fact, while in the second case, this is a nontrivial result obtained by D.~Romik and P.~\'Sniady and based on the papers~ \cite{KV}, \cite{VK}, and others. We comment on the proofs from the viewpoint of the theory of measurable partitions; another proof, using representation theory and generalized Schur--Weyl duality, will be presented elsewhere. We also study a new dynamics of Bernoulli variables on $P$-tableaux and find the limit  3D-shape of these tableaux.
\end{abstract}

\section{Introduction}\label{sec1}

This paper deals with common statistical and asymptotic properties of the classical Bernoulli scheme (i.e., a sequence of i.i.d.\ variables) and a generalized RSK (Robinson--Schensted--Knuth) algorithm. Simultaneously, we study the so-called combinatorial encoding of Bernoulli schemes. The combinatorial method of encoding continuous signals has certain advantages compared with the traditional invariant continual encoding; instead of functions with values in a continual space, it uses simple combinatorial relations between the coordinates of the message to be encoded. The main problem is to establish whether decoding is possible and how difficult it is.

The first example, originally considered in~\cite{V1}, is the encoding via Weyl simplices. It is typical, albeit quite simple; here we explain it in detail in terms of sequences of $\sigma$-algebras, or measurable partitions.

The second example, the encoding of Bernoulli schemes via the RSK correspondence, is much more complicated and interesting in itself, since it is related to the statistics of Young diagrams, asymptotic representation theory in general, etc. It appeared in the 1980s, in the paper~\cite{KV}, in connection with the study of central measures and, in particular, the Plancherel measure on Young tableaux. The main result said that the RSK correspondence, more exactly, its recording part (the $Q$-tableau) is a~homomorphism sending every one-sided Bernoulli measure  to some central measure on the set of paths in the Young graph. In particular, this homomorphism sends the one-sided Bernoulli shift to the so-called Sch\"utzenberger shift, which is an infinite counterpart of the ``jeu de taquin'' transformation introduced by Sch\"utzenberger.

In a remarkable recent series of papers \cite{RS}, \cite{S}, it was proved that the homomorphism considered in~\cite{KV} is an isomorphism, and thus, in terms of information theory, the encoding of Bernoulli schemes via RSK can be decoded.

We comment in detail on the principal line of~\cite{S} for the main (continuous) case, i.e., for the Plancherel measure. The decoding, i.e., the inverse isomorphism, is constructed in a complicated way, and, in our opinion, needs to be further analyzed. The authors of~\cite{S} use the limit shape theorems for Young diagrams proved in~\cite{KV},~\cite{LS} and identify, in a nontrivial way, this limit shape (more exactly, the curve~$\Omega$, see below) with the space of states of the Bernoulli scheme being decoded.

The important part here is not so much this identification and specific inversion formulas, as  the fact itself that the infinite $Q$-tableau contains all necessary information on a realization of the Bernoulli scheme. It is possible to prove this result by another method, relying more heavily on representation theory; it will be published separately.

There is another important characteristic of encoding: the growth of the number of states in the encoding sequence. In the first example, this growth is factorial,~$\{n!\}$, i.e., the number of states coincides with the order of the symmetric group~$S_n$; in the second example (the RSK correspondence), it coincides with the number of involutions in this group, i.e., is of order~$\{\sqrt {n!}\}$. It is not difficult to show that for a stationary combinatorial encoding, the growth of the number of states must be superexponential; however, the author does not know any meaningful examples of encodings with the growth of the number of states arbitrarily close to exponential.

Combinatorial encoding can also be applied to systems with a finite or countable set of states. Note that a related isomorphism  was considered in~\cite{VTs2}  from the viewpoint of representation theory  (the so-called concomitant representations).
In terms of the RSK algorithm, this case corresponds to discrete central measures and is much simpler than the case of the Plancherel measure. The decodability here is equivalent to the isomorphism, interesting in itself, between a Bernoulli shift and a Markov shift with a discrete set of states. A distinguishing feature of all these isomorphisms is that  the decoding has an infinite delay: to decode the first coordinate, we must know the coordinates of the encoding sequence with arbitrarily large indices. This effect is poorly studied in ergodic theory.

But what can then be said about the second component of the RSK algorithm, i.e., the $P$-tableau? It plays only an auxiliary part in the correspondence, and is not needed for decoding.

The main result of this paper (see Sect.~\ref{sec4}) shows that the $P$-tableaux have their own dynamics and asymptotics, which are related to deep and previously unknown properties of classical Bernoulli schemes. These properties are due to the existence of a linear order on the set of states. The RSK correspondence reveals the common
properties of the two orders on the states, this  linear order and the temporal order. Our theorem on the limit shape of $P$-tableaux shows an asymptotic interplay between these orders. It is somewhat unexpected that the limit shape itself coincides with the three-dimensional limit shape for diagrams, i.e., for $Q$-tableaux, which is  long known~\cite{Gr}. Apparently, the dynamics of $P$\nobreakdash-tableaux has not been studied earlier, and it raises some questions, which the author currently studies together with a group of young mathematicians. Some of the results are mentioned below.

We emphasize that this whole circle of problems is related to ergodic theory and information theory, that is why we discuss increasing sequences and decreasing sequences (i.e., filtrations) of partitions ($\sigma$-algebras), which inevitably appear in this context; their analysis explains the naturalness of results which seem  a priori mysterious.

Note also that all results on the asymptotic behavior of tableaux presented here have a direct interpretation as results on the asymptotic behavior of \textit{finite}
tableaux, and not only in the framework of Bernoulli schemes.

Combinatorial encoding leads to a number of new combinatorial and geometric problems on graphs and paths in graphs (see~\cite{V2}). The intriguing question here is what can be suggested as a counterpart of the RSK correspondence in the case of state spaces endowed with a partial order different from a linear one, or with some other structure. Almost nothing is known about this. But even for a linear order there are questions extending the domain of applicability of the RSK correspondence: limit shapes for other distributions of independent random variables, for Markov chains, etc.

In Sect.~\ref{sec2}, we consider the two examples of combinatorial encoding mentioned above; in Sect.~\ref{sec3}, we comment on the papers~\cite{RS}, \cite{S}, give an outline of the proof from these papers, and deduce combinatorial corollaries, in particular, for the representation theory of the infinite symmetric group. In Sect.~\ref{sec4}, we consider the dynamics and asymptotics of $P$-tableaux for a sequence of independent variables.

We emphasize again that there is a deep and poorly studied connection between the RSK correspondence and the theory of Schur--Weyl duality. It is this connection that is the source of various isomorphisms in the problems under study. It will be considered in a subsequent paper.

\section{Two fundamental examples}\label{sec2}

\subsection{Encoding of Bernoulli schemes via Weyl simplices}
We define several combinatorial objects on the infinite-dimensional cube $I^{\infty}\equiv I=\prod_{n=1}^{\infty} [0,1]$ endowed with the following product measure: $m^{\infty}\equiv m=\prod_{n=1}^{\infty} m^1$ where $m^1$ is the Lebesgue measure on the interval $[0,1]$. The one-sided shift  $T(\{x_k\})_n=x_{n+1}$, $n=1,2,\dots$, is an endomorphism of the measure space~$(I,m)$.

Consider two sequences   $\{\xi_n\}$ and $\{\eta_n\}$ of measurable partitions  of the space $(I,m)$.

An element of $\xi_n$ is a finite set, namely, an orbit of the group~$S_n$ acting by permutations of coordinates; thus, each element of~$\xi_n$ consists of all sequences $\{x_n\}\in I$ in which the coordinates with indices from~$n+1$ on coincide and the first $n$~coordinates coincide up to a permutation, i.e., have the same collection of elements. For convenience, we exclude from~$I$ the subset of zero measure consisting of all sequences that have at least two equal coordinates.

An element of $\eta_n$ is a set of finite measure (equal to $(n!)^{-1}$) containing all sequences
$\{x_n\}\in I$ in which the first $n$~coordinates are ordered in the same way and the coordinates with indices greater than~$n$ are arbitrary. It is natural to call such an element  a Weyl simplex, since its restriction to the $n$-dimensional subspace of the first coordinates is the intersection of an (open) Weyl chamber with the unit cube.

Recall that the space of classes of $\operatorname{mod}0$-coinciding measurable partitions of a~measure space (or the corresponding $\sigma$-subalgebras of the algebra of measurable sets) is a lattice with the following ordering\footnote{In combinatorics, the reverse order in the lattice of partitions of a finite set is used.}:

$\alpha \succ \beta$ if almost all elements of~$\alpha$ are composed (up to sets of zero conditional measure) of elements of~$\beta$.

The partition usually denoted by~$\epsilon$, the greatest element of the lattice, is the partition into singletons~$\operatorname{mod}0$; the zero (trivial, smallest) partition $\nu$ is the partition into a single nonempty set of full measure and the class of sets of zero measure.

The lattice of all measurable partitions is endowed with the topology induced by the weak operator topology in the $L^2$~space restricted to the operators of conditional expectation bijectively associated with $\operatorname{mod}0$ classes of  measurable partitions.

It is easy to verify the following properties of the partitions introduced above:

(a) the sequence of partitions $\{\xi_n\}_n$ is monotone decreasing, and the sequence~$\{\eta_n\}$ is monotone increasing\footnote{Monotone decreasing sequences of measurable partitions are sometimes called decreasing filtrations; if  a filtration converges to the trivial partition, then it is said to be ergodic.};

(b) for every  $n$, the partitions $\xi_n$  and $\eta_n$ are independent with respect to the measure~$m$;

(c) each of the two sequences of partitions is invariant under the one-sided shift~$T$:
$T\xi_n=\xi_{n+1}$, $T\eta_n=\eta_{n+1}$, $n=1,2,\dots$\,.

\begin{theorem}\label{th1}
The sequence $\{\xi_n\}_n$ converges, in the weak topology, to the trivial partition, i.e., is ergodic:
$\bigwedge_{n=1}^{\infty}\xi_n \equiv \xi_{\infty}=\nu$.

The sequence $\{\eta_n\}$ converges, in the weak topology, to the partition into singletons, i.e.,
$\bigvee_{n=1}^{\infty}\eta_n \equiv \eta_{\infty}=\epsilon$.
\end{theorem}

The first fact is well known: it means that the action of the infinite countable symmetric group on the infinite-dimensional cube by permutations of coordinates is ergodic.

Note that if an increasing sequence of measurable partitions $\{\alpha_n\}$ converges to the partition into singletons denoted by~$\epsilon$, then a decreasing filtration of measurable partitions $\{\beta_n\}$ such that each $\beta_n$ consists of independent complements to the elements of~$\alpha_n$, $n=1,2, \dots$, converges to the trivial partition~$\nu$.

The following fundamental fact of the theory of filtrations is less obvious: the converse is in general false, i.e., $\beta_n \rightarrow \nu$ does not imply $\alpha_n\rightarrow \epsilon$
    (see~\cite{VF}).

That is why, the second claim of the theorem, which says that almost all points of the Bernoulli scheme can be distinguished by the list of pairwise inequalities, does not follow from the first claim, the ergodicity of the action of~$S_{\infty}$. But it is not difficult to prove it directly (see~ \cite{V1}). There are various combinatorial applications of this fact, for example, the following paradoxical interpretation: almost all infinite Weyl simplices consist of a single point, since the list of all inequalities between all coordinates uniquely determines the simplex and, according to the theorem, determines almost every sequence\footnote{And almost every, with respect to the standard Gaussian measure, Weyl chamber consists of a~single ray. Here is a similar paradox: the partition of the infinite-dimensional space with a~Gaussian measure into the rays starting at the origin is the partition into singletons.}. Of course, it is important here that we consider a set of measure~$1$.

Consider the space $W=\prod_{n=1}^{\infty} \mathbf{n}$, $\mathbf{n}=\{1,\dots,n\}$, i.e., the space of all sequences~$\{z_n\}$ where $z_n\in \mathbf{n}$, and denote by~$\theta$ the Borel measure on~$W$ equal to the product of the uniform measures on the factors $\mathbf{n}$. Note that $W$~should be regarded as the space of all paths in the graph $\Gamma=\bigcup_{n=1}^{\infty} \mathbf{n}$ where any two neighboring levels, i.e., $\mathbf{n}$ and $\mathbf{n+1}$, form a complete graph. In~\cite{V1}, the space~$W$ was called the triangular compactum.\footnote{This compactum, under the name of ``the space of virtual permutations of positive integers,'' appeared in another connection, namely, it was defined in~\cite{KOV} as the inverse spectrum of the sequence of symmetric groups~$\{S_n\}$ with respect to the projection forgetting the last coordinate of a~permutation.}

\textit{The combinatorial meaning of the above example of encoding and Theorem~{\rm\ref{th1}} is as follows: there is a nontrivial isomorphism~$\tau$ between two dynamical systems  $(I,m,T)$ and $(W,\theta, S)$},
where
$T$~is the one-sided shift on the space~$I$ (a~Bernoulli endomorphism). The map $\tau\colon I\rightarrow W$, $\tau(\{x_n\})=\{z_n\}$, sends a~sequence $\{x_n\}\in I$ to the sequence
$\{z_n\}\in W$ where $z_n=\#\{i:1\leq i\leq n,\,x_i\leq x_n\}$, and Theorem~\ref{th1}  says that $\tau$ is an isomorphism $\operatorname{mod}0$ between the continual and the discrete dynamical systems (a shift of continual objects and a shift of combinatorial objects).

The shift $S$ in the space $W$ is defined by the formula $S=\tau\cdot T\cdot \tau^{-1}$ and is a~shift in the path space of a graph (for details, see~\cite{V1}). In more detail: ${S(\{z_n\})=\{z'_n\}}$ where $z'_n=z_{n+1}$ if $x_1>x_{n+1}$ and $z'_n=z_{n+1}-1$ if $x_1<x_{n+1}$. In turn, the condition ${x_1 \gtrless x_{n+1}}$ can be expressed in terms of $z_1,\dots,z_n,z_{n+1}$ or in terms of the corresponding permutation. The important thing is that the functions $d_n$, ${n=1,2,\dots}$, on~$I$, where
$$
d_n(\{x_n\})=\#\{i<n:x_1>x_i\},
$$
are well defined as functions on~$W$, i.e., as functions of $\{z_n\}$, and, in particular, a.e.\ on
$(W,\theta)$ there is a limit
$$
\lim_{n \to \infty}\frac{d_n}{n},
$$
equal to the density of the set coordinates of $\{x_n\}$ less than~$x_1$, which in our case is exactly~$x_1$.

It is this fact that allows us to find the inverse isomorphism  $\tau^{-1}\colon W\rightarrow I$, which can be easily written explicitly if we first recover the space of states of the Bernoulli system corresponding to the first (and hence any) coordinate. Here, this space is the interval
$[0,1]$, and it arises as the collection of all limits of the quotients~$d_n/n$, i.e., as the ``limit shape'' of states of the discrete system
$(W,\theta,S)$. ``Time reversal'' occurs: the space of initial states of the Bernoulli scheme arises as a limiting space.

The map $\tau$ is defined as a homomorphism from the triple~$(I,m,T)$ onto the quotient space~$I/{\eta_{\infty}}$, but,  in view of the equality $\eta_{\infty}=\epsilon$ from Theorem~\ref{th1}, it turns out to be an isomorphism.  However, the formula for $d_n$ and the existence of $\lim_{n\to\infty}d_n/n$ provide a direct proof of the inversion formula for~$\tau$.

The above construction  relies on the decompositions, consistent for different~$n$, of the $n$-dimensional cubes with the Lebesgue measure into the direct product $\operatorname{mod}0$ of two spaces: a typical Weyl simplex~$\Delta(n)$ (for example, the simplex of monotone sequences) and the symmetric group~$S_n$ (the Weyl group). This allowed us to construct two shift-invariant sequences of partitions of the infinite-dimensional cube, one increasing and the other decreasing, such that for every $n$ they are independent complements of each other, the increasing sequence converges to the partition into singletons, and the decreasing one converges to the trivial partition. One can suggest other similar examples, for instance, constructed from other series of simple Lie groups.

However, a more important thing is the relation to the theory of filtrations.

\begin{lemma}\label{lm1}
The  filtration consisting of the decreasing sequence of orbits of the groups~$S_n$, $n=1,2,\dots$, is standard in the sense of the theory of filtrations (see~\cite{VF}), i.e., isomorphic to the $n$-adic filtration of $\prod_{n=1}^{\infty}\mathbf{n}$.
\end{lemma}

This lemma follows from Theorem~\ref{th1}; in fact, an isomorphism between the cube~$I$ and the triangular compactum~$W$ defines a system of complementary partitions needed to establish the standardness. We emphasize that the decodability is equivalent to the standardness of some filtration.

That is why, the constructed nontrivial isomorphism between the actions of the groups $S_{\infty}$ and $\sum_{n=1}^{\infty} \mathbb{Z}/n$ on the cube $I$, which is orbital for every~$n$, should serve as a model for other situations. In the next section we consider such a~situation.

\subsection{Encoding via the RSK correspondence}
Here we generalize the previous construction using the paper~\cite{KV} and the subsequent papers~\cite{RS}, \cite{S}.

Much wider opportunities open if we allow a more arbitrary construction of partitions of the finite-dimensional cubes~$I^n$ used to define two sequences, one decreasing and the other increasing, of measurable partitions of the infinite-dimensional cube, the first one converging to the trivial partition and the second one, to the partition into singletons; here, the $n$th elements of the sequences  must not necessarily be  independent complements of each other on the whole cube with respect to the Lebesgue measure for every~$n$.

In the next example, an important role is played again by the symmetric groups~$S_n$, $n=1,2,\dots$, and the well-known Robinson--Schensted--Knuth  (RSK) correspondence between~$S_n$ and the pairs of Young tableaux with the same diagram. We assume this correspondence known (see, e.g., \cite{St,Kn}), but consider it, following~\cite{KV},  not only for permutations, but for sequences of elements of an arbitrary linearly ordered set, in particular, the interval~$(0,1)$, i.e., for elements of the infinite-dimensional unit cube~$I$ endowed with the Lebesgue measure~$m$.

\textit{Apply the} RSK \textit{correspondence first to a vector $x^n=(x_1,\dots,x_n)$ of the $n$-di\-men\-sion\-al cube
 $[0,1]^n=I^n$. This gives two Young tableaux of the same shape, the insertion tableau~$P(x^n)$ and the recording tableau~$Q(x^n)$}; the cells of the first tableau are filled with the numbers~$x_i$, $i=1,\dots, n$, inserted according to a certain rule, while the second tableau is a \textit{standard Young tableau} (its cells are filled with the positive integers from~$1$ to~$n$; this tableau records the process of growth of the diagrams). In contrast to the case of~$S_n$ (in which vectors~$x^n$ are permutations),  the $P$- and $Q$-tableaux here are very distinct in nature and play substantially different roles in what follows. Below we will analyze the quite different asymptotic behavior of these tableaux.

In order to extend the RSK correspondence to infinite sequences, we must keep track of the asymptotic behavior of the correspondence for finite vectors. For this, we must study various partitions of the cube~$I$ related to the RSK correspondence.

Denote the set of all Young diagrams with $n$~cells by~$\Lambda_n$; with each diagram ${\lambda\in \Lambda_n}$, we  associate the subset~$I(\lambda)$ of finite measure in the $n$-dimensional cube~$I^n$ consisting of all vectors for which the (coinciding) shape of the $P$- and $Q$-tableaux coincides with~$\lambda$. Keeping the same notation, we regard~$I(\lambda)$ as cylinder sets in the infinite-dimensional cube~$I$. Thus, for every $n$ we obtain a finite partition of the cube~$I$, with the number of elements equal to Euler's partition function~$p(n)$; denote it by~$\rho_n$. Obviously, $\rho_n \prec \eta_n$, since every element of~$\rho_n$ is the union of some Weyl simplices, namely, the simplices whose ordering corresponds to the given diagram~$\lambda$.
The sequence of partitions $\{\rho_n\}$ is neither decreasing nor increasing. The interrelations between the elements of these partitions for $n$~and~$n+1$ and on the whole are described precisely by the Young graph.

More important partitions are the products $\bigvee_{k=1}^n \rho_k \equiv \bar\rho_n $.
By definition, the sequence $\{\bar\rho_n\}$ is increasing, and every its element is the union of all vectors $\{x_i\}_i$ from the cube $I$ that have the same $Q$-tableau of the initial segment  $x^n\equiv(x_1,\dots,x_n)$. In other words, the elements of~$\bar\rho_n $ are indexed by the Young tableaux with $n$~cells. \textit{The encoding of the Bernoulli scheme via the} RSK \textit{correspondence sends a realization~$\{x_n\}_n$ of the Bernoulli scheme to the sequence of Young tableaux $Q(x^n)$, or, since this sequence stabilizes, to an infinite Young tableau.} The distinguishability (decodability) problem is whether or not $\bigvee_{k=1}^\infty\rho_k \equiv \bar\rho_\infty$ is the partition into singletons. This is essentially the subject of the papers \cite{RS}, \cite{S}, as well as the next section. In these papers, it is established that $\bar\rho_\infty$ is the partition into singletons $\operatorname{mod}0$. In other words, decoding is possible (distinguishability holds)\footnote{In \cite{S}, the term ``asymptotic determinism'' is used, which is less natural from the viewpoint of signal transmission.}. As to
$\lim_n\xi_n$, it follows from Theorem~\ref{th3} and the remark on the relation between increasing and decreasing sequences of mutually complementary partitions that this limit coincides with the trivial partition~$\nu$. In the next section (Sect.~\ref{subs3.3}), we compare these partitions in more detail.

\section{Statements and corollaries of the main theorems}\label{sec3}

In this section, we comment on the papers \cite{KV}--\cite{S} and deduce corollaries from the main theorems.

\subsection{Statements}
The first result, obtained in \cite{KV}, connects central measures on the Young graph with Bernoulli schemes:

\begin{theorem}[\cite{KV}]\label{th2}
The quotient of the infinite-dimensional cube $(I,m)$ endowed with the product measure~$m$ by the limit
$\bigvee_{n=1}^{\infty}\bar\eta_n$  of the increasing sequence of partitions $\{\bar\eta_n\}$
is the space~$\mathcal{T}$ of all infinite standard Young tableaux endowed with the Plancherel measure~ $\mu$.
\end{theorem}

It is easy to conclude that the image of the one-sided shift~$T$ on the cube~$I$ under the above identification is an infinite version of the Sch\"utzenberger transformation~Sch (see~\cite{St}). The proof in~\cite{KV} consisted essentially in directly  verifying that the projections of product measures are central measures. In itself, this did not provide a new proof of the list of central measures (Thoma's theorem).

A very important result obtained in a series of recent papers~ \cite{RS},~\cite{S} significantly refines the previous one and consists in the following:

\begin{theorem}[\cite{RS},~\cite{S}]\label{th3}
The homomorphism  $\pi\colon(I,m,T)\rightarrow(\mathcal{T},\mu,\operatorname{Sch})$ is an isomorphism. In other words,  a sequence of independent random variables can be recovered (decoded) from an infinite Young tableau.
\end{theorem}

Below, we outline the proof given in \cite{RS},~\cite{S} and trace an analogy with our first example. We will speak about the Plancherel measure only; a simpler analysis of other central measures on the space of Young tableaux follows the same plan and will be considered later.

The established isomorphism leads to useful corollaries for the theory of representations of the infinite symmetric group and, in a sense, settles the question of the asymptotic behavior of the recording $Q$-tableaux in the infinite RSK algorithm (see the subsequent sections).

The asymptotic behavior of the insertion $P$-tableaux has not apparently been studied so far; we analyze it in Sect.~\ref{sec4}.

\subsection({The idea of the proof (\000\134cite\{RS\}, \000\134cite\{S\}): two descents}){The idea of the proof (\cite{RS}, \cite{S}): two descents}
As mentioned above, Theorem~\ref{th2} is proved in~\cite{KV} by directly calculating the measures of cylinders with respect to the image of the Lebesgue measure under the projection~$\pi$. This implies a~homomorphism of dynamical systems between the one-sided Bernoulli shift with a continual set of states and the Sch\"utzenberger shift with respect to the Plancherel measure on the space of infinite Young tableaux.

Note that this fact alone already implies that the Sch\"utzenberger shift is \textit{isomorphic} to the one-sided Bernoulli shift: this follows from the general theorem saying that a quotient of a Bernoulli shift having the same entropy is isomorphic to the original shift. However, Theorem~\ref{th3} provides a natural specific isomorphism~$\pi$.

Let us explain the essence of the proof of this theorem from~ \cite{RS},~\cite{S}, as  understood by the author of this paper. To establish the isomorphism, one must prove the invertibility of the homomorphism~$\pi$ from the space~$(I,m)$ onto the space $(\mathcal{T}, \mu)$ of infinite Young tableaux endowed with the Plancherel measure~$\mu$. In simpler words, one must prove that, for instance, the first element of the sequence $x=\{x_1,x_2,\dots\}$, i.e., $x_1$, can be recovered from its image under the homomorphism~$\pi$, i.e., from an infinite Young tableau, more exactly, from the recording $Q$-tableau.

Actually, this theorem states the existence of a measurable map inverse to the projection~$\pi$, i.e., the decodability of a sequence $\{x_n\}$ from the sequence of its recording tableaux (or from the corresponding infinite Young tableau). Since the construction is shift-invariant, having recovered the first element, in the same way one can recover the second one, etc., that is, decode a realization of the Bernoulli scheme from its encoding via the RSK correspondence.

In contrast to our first example (encoding via Weyl simplices, see above), where the decoding is given by the limiting inversion formula $x_1=\lim_n d_n(z)/n$, in this case the decoding (inversion formula) presented in
\cite{RS}, \cite{S} is much more involved, though the general context resembles the simplest case. I will explain not so much the formula itself (it is not very useful and important), as the very existence of inversion, i.e., decoding. A nontrivial construction of such an inverse map is precisely the content of the papers~\cite{RS},~\cite{S}.

In the RSK algorithm, the first element~$x_1$ of the sequence behaves in the $P$-tableau as follows: at the first step, it is placed into the first cell of the first column, i.e., the cell~$(1,1)$, and then it can move only down the first column.

We know from the limit shape theorem that the length of the first column (and the first row) grows as $2\sqrt{n}$ where  $n$~is the number of cells. As we will see, the following result holds:

\begin{lemma}\label{lm2}
The limiting position of the first element $x_1=\bar x$ in the first column as $n \to \infty$ is $\bar x\cdot 2\sqrt{n}$; in other words, the quotient of the number of the row containing~$x_1$ after $n$~steps of the RSK algorithm to the length of the first column tends to the value~$x_1$ itself.
\end{lemma}

Therefore, one can recover the value~$x_1$ by finding the normalized number of the row containing it. But one is allowed to use only the $Q$-tableau. It turns out that this can be done, which is the key point of the proof; it relies on finer properties of the infinite RSK algorithm, namely, on information about the numbers occupying certain parts of the $P$- and $Q$-tableaux.

The point is that the $Q$-tableau has its own descent dynamics related to the Sch\"utzenberger transformation, and this descent proceeds with the same velocity. Therefore, by calculating this velocity we can obtain the desired value~$x_1$.

\begin{figure}[h]
\begin{center}
\includegraphics[width=0.35\textwidth]{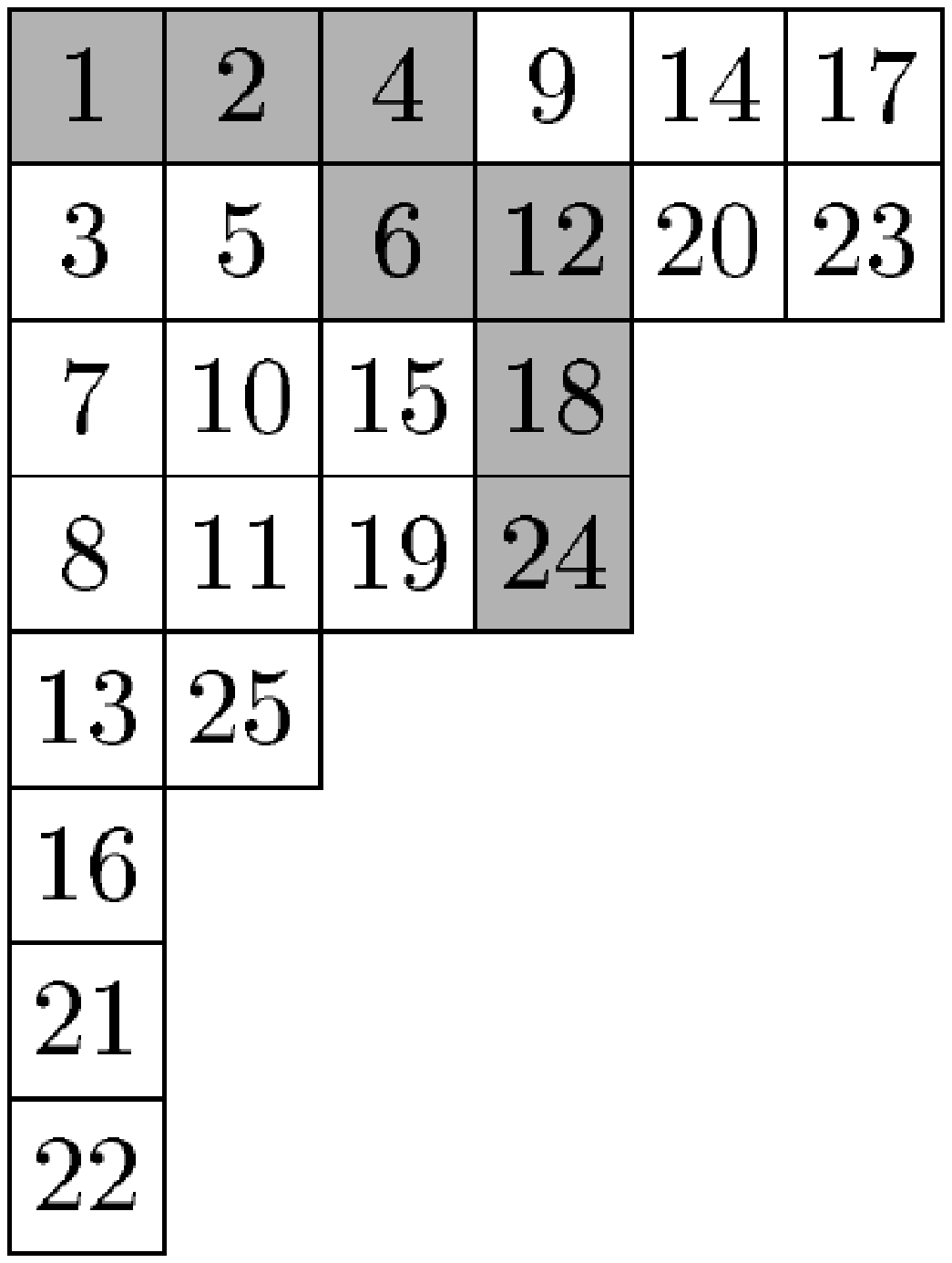}
\end{center}
\caption{The nerve of a tableau (numeration)}
\label{fig1}
\end{figure}

The descent in the $Q$-tableau is defined via the Sch\"utzenberger shift as follows. An infinite Young tableau should be regarded as a numeration of the lattice~$\mathbb{Z}_+^2$ by the positive integers (see~\cite{VB}). Almost every, with respect to the Plancherel measure, Young tableau fills the whole lattice. Assume that a Young tableau is fixed. Consider the connected chain of cells  in~$\mathbb{Z}_+^2$ starting from the cell~$(1,1)$ and constructed according to the following rule: if we have already constructed the chain up to some cell, then
the next cell  is either its right neighbor or its bottom neighbor, namely, whichever of the two contains the smaller number. Thus, every infinite tableau gives rise to a connected sequence of lattice cells, which, together with the numbers occupying these cells, will be called the \textit{nerve of the Young tableau} (see Fig.~\ref{fig1}).

The Sch\"utzenberger shift sends a Young tableau~$t$ to the tableau obtained from~$t$ as follows: fill the first cell $(1,1)$ of the nerve of~$t$ with the number occupying the second cell of the nerve (which is always $2$), fill this second cell  with the number occupying the third cell of the nerve, etc.; the numbers in all cells outside the nerve remain the same as in~$t$. Now decrease the numbers in all cells by~$1$. Obviously, we obtain a new Young tableau. It is the image of the original tableau~$t$ under the Sch\"utzenberger shift~$\operatorname{Sch}$.

Let $t$ be an infinite Young tableau. Consider its finite fragment~$t_n$ consisting of $n$~cells; denote the lattice coordinate of the last cell of the nerve of~$t_n$ by $(a_1(t,n),a_2(t,n))$.
It follows from the known facts on the asymptotic behavior of tableaux with respect to the Plancherel measure that the growth of these coordinates  is of order~$\sqrt{n}$. The main result of the papers~\cite{RS},~\cite{S} is essentially as follows:

\begin{lemma}\label{lm3}
For almost all, with respect to the measure~$m$, realizations $\{x_1,x_2, \dots\}$, the first coordinate can be expressed as
$x_1=\lim_n a_1(t,n)/\sqrt{n}$.
\end{lemma}

Taking into account Lemma~\ref{lm1}, we may say that the growth rate of the number of the row in the $P$-tableau containing the first element~$x_1$ coincides with the growth rate of the first coordinate of the cell of the nerve in the $Q$\nobreakdash-tableau.

On the other hand, from the general theory of the RSK correspondence (see~\cite{St}) we know what numbers occupy the first row of the $P$-tableau; that is why, it is useful to transpose the tableau (which corresponds to inverting the order of the sequence) and keep track of the first row, which is the image of the first column. This leads to cumbersome and possibly unnecessary calculations. In the author's opinion, the problem of finding a simple proof of this lemma is still open. The authors of~\cite{RS},~\cite{S} use the quotient $a(t,n)/b(t,n)$ of the coordinates of the end of the nerve, prove that this quotient has a limit, and then
identify the limit of the nerve with the point of the limit shape of Young diagrams (\cite{VK}, \cite{LS}) lying on the ray corresponding to this quotient of coordinates. In this way, the space of states of the Bernoulli scheme gets identified with the limit curve~$\Omega$. This argument heavily relies on the theorem about the limit shape of Young diagrams with respect to the Plancherel measure (\cite{VK}, \cite{LS}), which does not seem indispensable. Above, we kept the interval as the space of states, which, clearly, is not important. Of course, the relation between the dynamics of the $P$- and $Q$-tableaux is not simple, and it is at the heart of the problem.

Another plan of proving the isomorphism relies on directly establishing an equivalence of two representations of the infinite symmetric group, or, equivalently, establishing an isomorphism of shifts on measure spaces, one with a~central measure, and the other one with a Bernoulli measure.

\subsection{Combinatorial corollaries}\label{subs3.3}
In this section, we present some combinatorial corollaries and interrelations between the partitions considered above.

Here is the main corollary, which allows one to compare both examples.

\begin{corollary}{\rm(distinguishability of sequences in terms of the parameters of the diagrams)}
Almost every realization $\{x_n\}_n$ is uniquely determined by the parameters of the Young diagrams corresponding to the tableaux~$\{t_i\}$ obtained as the $Q$-tableaux in the  \textup{RSK} correspondence, for instance, by the collections of all row lengths, or all column lengths, or numerations of the lattice~$\mathbb{Z}_+^2$, etc.
\end{corollary}

In other words, from every family of these collections, for almost every tableau with respect to the Plancherel measure, one can recover the unique realization of the Bernoulli scheme generating this tableau under the RSK correspondence. The above collections are included as a small part into the collection of all pairwise inequalities between coordinates, which we discussed in Sect.~2.

Clearly, at a finite level there is no such thing: one cannot recover a finite permutation from its $Q$-tableau, though it can be recovered from the list of all inequalities. It would be very interesting to find other such encodings of Bernoulli schemes.

Consider the partition of the finite-dimensional cube~$I^n$ into the sets of vectors with the same $P$-tableau (insertion tableau), and carry it over to the cube, obtaining a partition~$\bar\xi_n$, as follows:

\begin{definition}
An element of $\bar\xi_n$ consists of all vectors $\{x_n\}\in\nobreak I$ in which the coordinates with indices greater than~$n$ coincide and the initial fragments~$x^n$ have the same $P$-tableau.

Recall that an element of the partition~$\bar\rho_n$ consists of all vectors having the same $Q$-tableau  (recording tableau).
\end{definition}

Let us compare these partitions with the Weyl partitions~$\xi_n$, $\eta_n$ (see Sect.~\ref{sec2}).

Of course, $\bar\rho_n\prec \eta_n$, since an element of~$\bar\rho_n$ is composed of Weyl simplices. On the other hand, one can easily deduce from the definition of $P$-tableaux that elements of~$\bar\xi_n$ refine elements of~$\xi_n$, i.e., $\bar\xi_n \succ \xi_n$. The latter is the partition into the orbits of the symmetric group.

As observed in Sect.~\ref{sec2}, the sequences of partitions~$\{\eta_n\}_n$ and~$\{\xi_n\}_n$ are termwise (i.e., for every~$n$) mutually independent, the first one converges to the partition into singletons, and the second one, to the trivial partition.

It follows from Theorem~\ref{th3} and some additional considerations that, though we have coarsened~$\{\eta_n\}_n$
and refined~$\{\xi_n\}_n$, and though $\bar\rho_n$~and~$\bar\xi_n$ are independent only when restricted to the elements of the partition~$\rho_n$, the following result still holds.

\begin{proposition}\label{pr1}
We have
$$
\lim_n \bar\rho_n =\epsilon, \qquad \lim_n \bar\xi_n=\nu.
$$
\end{proposition}

In contrast to the first example, the partitions $\bar\xi_n$~and~$\bar\rho_n$ are not independent complements of each other; however, as follows from definitions, their restrictions to the elements of the partition~$\rho_n$ are. Indeed, for a finite~$n$ and a fixed diagram, the $P$-~and~$Q$-tableaux are, obviously, independent with respect to the uniform measure on the set of permutations with the given diagram. This implies

\begin{corollary}
The partition $\bar\xi_n$ of the cube~$I$ is a subpartition of the partition into the orbits of the action of the group~$S_n$ by permutations of coordinates and is an independent complement to the restriction of the partition~$\bar\rho_n $ to the elements of the partition~$\rho_n$.
\end{corollary}

The answer to the question of whether there exists a natural decreasing sequence of independent complements to the sequence~$\{\bar\rho_n\}_n $ seems to be negative; however, the above partial independence suffices to deduce that the first equality in Proposition~\ref{pr1} implies the second one.

Let us link these partitions to combinatorial definitions of equivalence of permutations.

Recall that two finite sequences $\omega_1$, $\omega_2$ of elements of an arbitrary linearly ordered set $L=\{a,b,c,\dots\}$ (for example, two permutations of the set $L=\{1,\dots,n\}$) are \textit{Knuth equivalent} if each of them can be obtained from the other one by an arbitrary sequence of applications of the following operations on triples of neighboring coordinates: if
$a<b<c$, then
$$
bac \sim bca\quad\text{and}\quad acb \sim cab.
$$

Clearly, these operations do not change the collection of coordinates.

Also recall (see \cite{St}, \cite{Fo}) that two permutations of elements of a linearly ordered set are Knuth equivalent if and only if their insertion tableaux (i.e., $P$-tableaux) in the RSK correspondence coincide. The Knuth equivalence is, of course, stronger that the orbit equivalence of the action of the symmetric group by permutations of coordinates. There is also a notion of the dual Knuth equivalence on~$S_n$, which is defined as the Knuth equivalence of the inverse permutations (see~\cite{Fo}). In other words, permutations are dual Knuth equivalent if and only if their recording tableaux ($Q$-tableaux) coincide. Both notions can be applied to vectors of the finite-dimensional cube~$I^n$, meaning that the dual equivalence identifies vectors ordered in the same way.

\begin{proposition}\label{pr2}
An element of the partition $\bar\xi_n$ is a class of sequences  whose coordinates coincide for indices greater than~$n$ and  fragments~$x^n$ are Knuth equivalent.

An element of the partition $\bar\rho_n$ is a class of sequences whose initial fragments are dual Knuth equivalent.
\end{proposition}

In our framework, the Knuth equivalence is ergodic, while the dual Knuth equivalence distinguishes realizations. The second claim of Proposition~\ref{pr2} is equivalent to Theorem~\ref{th3}.

Essentially, the first question is about the convergence of the empirical distribution of the first coordinate of the initial fragment of the sequence averaged over the Knuth equivalence to the uniform distribution. Averaging over the whole symmetric group gives the uniform distribution by the martingale convergence theorem, or simply by ergodicity. But here, instead of an orbit of the symmetric group, we have the much narrower Knuth equivalence relation. The corresponding filtration, as we have seen, is ergodic and standard, which follows from Theorem~\ref{th3}. By the same theorem, the sequence~$\bar\rho_n$ converges to the partition~$\epsilon$ into singletons and corresponds to the dual Knuth equivalence (the equality of recording tableaux). It is more poorly understood\footnote{Computer simulations, performed by my student P.~Naryshkin, show that the convergence of the Knuth equivalence to the trivial partition is very slow. It would be interesting to describe a~natural countable group for which the orbit partition coincides with the Knuth equivalence.}.

There remain many interesting questions about infinite and asymptotic counterparts of combinatorial facts related to finite partitions and equivalences in the theory of symmetric groups and Young tableaux, such as the Knuth equivalences, Sch\"utzenberger involution, Littlewood--Richardson rule, etc. It seems that the difference between probabilistic and finite statements will be very large, however, ideological parallels are extremely important. We conclude this section by observing that the analysis of two sequences of partitions corresponding to the (direct and dual) Knuth equivalences makes quite natural the fact that in the limit of finite RSK correspondences with respect to just one component (the $Q$-tableaux) decoding is possible.

\subsection{The basic representation of the infinite symmetric group and semidirect products}
Beside corollaries related to encoding of Bernoulli schemes, there is another interesting corollary of the isomorphism between the space of infinite Young tableaux with the Plancherel measure and the infinite-dimensional cube~$(I,m)$ with the Lebesgue measure, namely, a new realization of an important irreducible representation, called the
\textit{basic} representation, of the infinite symmetric group in the Hilbert space $L^2(I,m)$ of square integrable functions on~$I$.

Recall that the $C^*$-group algebra  $C^*[S_{\infty}]$ of the infinite symmetric group has a~natural structure of an $AF$-algebra, i.e., is the cross product of a commutative algebra, which is the Gelfand--Tsetlin algebra of bounded measurable functions on the space of infinite Young tableaux (this space is the spectrum of the Gelfand--Tsetlin algebra), and the so-called tail equivalence relation on the space of infinite Young tableaux. More exactly, the group algebra has a~groupoid realization. That is why, the regular representation of the group~$S_{\infty}$ is a hyperfinite $\operatorname{II}_1$ factor, which has a~well-known von Neumann realization.

The tail equivalence relation can be described in a more familiar way as the orbit equivalence relation for certain groups, e.g., for the adic automorphism (Young automorphism, see~\cite{Vad},~\cite{Vad1}) or a group generated by involutions~\cite{VTs1}. But every factor representation of a cross product has a diagonal (Koopman) irreducible analog. Usually, it is a Koopman representation that is defined first, and then a~von Neumann factor is constructed from it; but here we move in the opposite direction. In our case, this analog is constructed as a representation of the group algebra~$C^*[S_{\infty}]$ in the Hilbert space of functions on the space of infinite Young tableaux that are square integrable with respect to the Plancherel measure. The resulting representation is irreducible, because the tail equivalence relation is ergodic and the weak closure of the Gelfand--Tsetlin algebra is a~maximal commutative subalgebra in the algebra of operators. All these facts are obvious and well-known. This Koopman representation of the cross product will be called the \textit{basic}, or \textit{Plancherel, representation} of the infinite symmetric group.\footnote{In a similar way,  Koopman representations can be constructed for an arbitrary central measure of any $AF$-algebra and, in particular, the group algebra of the infinite symmetric group. Such representations are called  \textit{concomitant} (see~\cite{VTs2}).}

The above isomorphism allows us to realize the basic representation of the infinite symmetric group in the space~$L^2(I,m)$. In this model, it is not a~permutation representation, i.e., is not generated by any action of the symmetric group on the cube~$(I,m)$ itself. Explicit formulas for the action of the Coxeter generators can be obtained  via the RSK correspondence.

Clearly, in this model the basic representation is the weak limit of a sequence of reducible representations of the groups~$S_n$ in the spaces~$L^2(I^n,m^n)$, and one may try to prove the irreducibility of the limit representation directly. This approach is closely related to generalizations of the Schur--Weyl duality to infinite-dimensional groups (see~\cite{VTs}). These questions will be considered in more detail in another publication.

\section{Bernoulli schemes and Young tableaux}\label{sec4}

We turn to the second component of the RSK correspondence, the sequence of insertion tableaux, which, at first sight, is left out of the previous considerations, since the second (recording) component alone allows one to invert the infinite RSK algorithm. Moreover, an explanation is needed for the fact that in the finite version of the RSK algorithm, both tableaux have equal rights and, of course, a~permutation cannot be recovered from only one of the tableaux, $P$~or~$Q$, alone. In other words, the inversion theorem is distinctively asymptotic.

\subsection({The movement of coordinates in \$P\$-tableaux: from the first row to the first column}){The movement of coordinates in $P$--tableaux: from the first row to the first column}
In contrast to the $Q$-tableaux, the $P$-tableaux have no nontrivial limit  in our scheme; more exactly, the following property holds:

\begin{proposition}\label{pr3}
For almost all, with respect to the measure~$m$, realizations\break $\{x_n\}_n\in I$, the sequence of values inserted into the cell~$(1,1)$ of the $P$-tableau tends to zero as
$n\to \infty$; the same holds for every cell. In other words, the sequence of $P$-tableaux weakly converges to the zero tableau for $m$-almost every realization.
\end{proposition}

\begin{pf}{Proof}
The number occupying the cell $(1,1)$ of the $P$-tableau after $n$ steps of the RSK algorithm is $\min_{1\leq i \leq n}x_i$, which tends to zero.  For every cell~$(i,j)$, the sequence of coordinates inserted into it is monotone decreasing, and if the limit of this sequence were positive, this would mean that the number in every cell of the infinite rectangle $(i+q,j+p)$, $q,p=1,2,\dots$, is strictly greater than this limit; but the set of indices of the coordinates inserted into the cells of this rectangle has density~$1$, and for every subsequence of indices of density~$1$,
the corresponding values are uniformly distributed on the interval.
\qed\end{pf}

\begin{proposition}\label{pr4}
For almost every realization $\{x_1,x_2,\dots \}$ and for every its coordinate $x_m$, $m\geq 1$, there exists a number $M\gg m$ such that the $P$-tableau of the fragment $\{x_1,\dots,x_{m+M}\}$  contains~$x_m$
in the first column.
\end{proposition}

\begin{pf}{Proof}
Recall that  after every insertion, the element occupying a given cell of the $P$-tableau either remains in the same cell or is bumped into a cell of the next row to the left of the original cell, i.e., the number of the column containing this element does not increase. Besides, this process of bumping is infinite, because there are coordinates with values arbitrarily close to the given one. A careful analysis shows that an element cannot remain in the same column with number greater than~$1$ infinitely often. In other words, in finitely many steps it reaches the first column.
\qed\end{pf}

My students I.~Azangulov and G.~Ovechkin have recently refined this result (see~\cite{AO}).

\begin{proposition}\label{pr5}
For any $\epsilon>0$ and $n\in\Bbb N$ there exists a set of realizations with measure greater than
$1-\epsilon$ such that the $n$th coordinate of every realization from this set reaches the first column of the $P$-tableau  after  $c(\epsilon)n^2$ steps of the RSK algorithm; here $c(\bcdot)$ is a positive function.
\end{proposition}

Note that this upper bound is close to being sharp; indeed, as proved in~\cite{VK}, the asymptotic length of the first column (and the first row) is equal to $2\sqrt n$, so the exponent~$2$ cannot be decreased.

\begin{figure}[h]
\begin{center}
\includegraphics[width=0.5\textwidth]{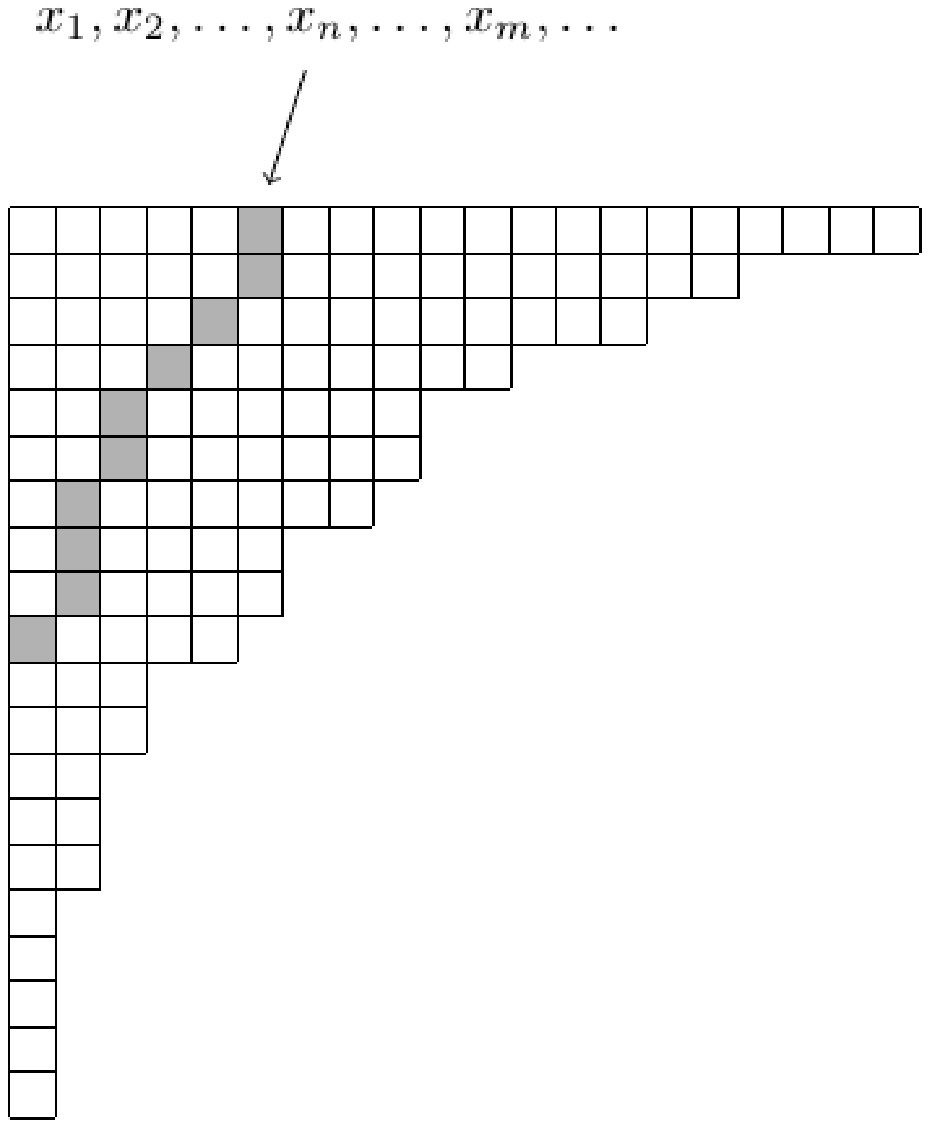}
\end{center}
\caption{The path of a random coordinate in the $P$-tableau}
\label{fig2}
\end{figure}

Thus,  on the $P$-tableau, every coordinate of almost every realization moves from a cell of the first row to a cell of the first column in time of order~$n^2$ (see Fig.~\ref{fig2}). Afterwards, it continues moving down  the first column.

This movement of coordinates  along the $P$-tableau requires further study: one must estimate the velocity of the first coordinate along the first column (this is important for the proof of the inversion theorem, see Theorem~\ref{th3}).  Note that it follows from above that the numbers occupying the first column
at every moment can be regarded as the support of the empirical distribution of some finite fragment (with volume of order~$\sqrt n$) of a given realization, which converges to the Lebesgue measure (i.e., to the theoretical distribution).

Thus, the $P$-tableau changes with~$n$, but, unlike the $Q$-tableau, does not stabilize, and its dynamics is not exhausted by the  movement described above.

\subsection({The evolution of \$P\$-tableaux, the limit shape, and the arch}){The evolution of $P$-tableaux, the limit shape, and the arch}
Consider the question of how the number occupying a given cell of the $P$\nobreakdash-tableau changes. For example, consider the behavior of the number occupying the last cell of the first row or the first column of the $P$-tableau. It is clear that the sequence of these numbers converge to~$1$. In the same way, we can take a~cell $(i\sqrt n, j\sqrt n)$ and consider the sequence of  values inserted into this cell and the limit of these values. It turns out that it is more convenient to use polar coordinates for the tableau, see below.

\begin{figure}[h]
\begin{center}
\includegraphics[width=0.6\textwidth]{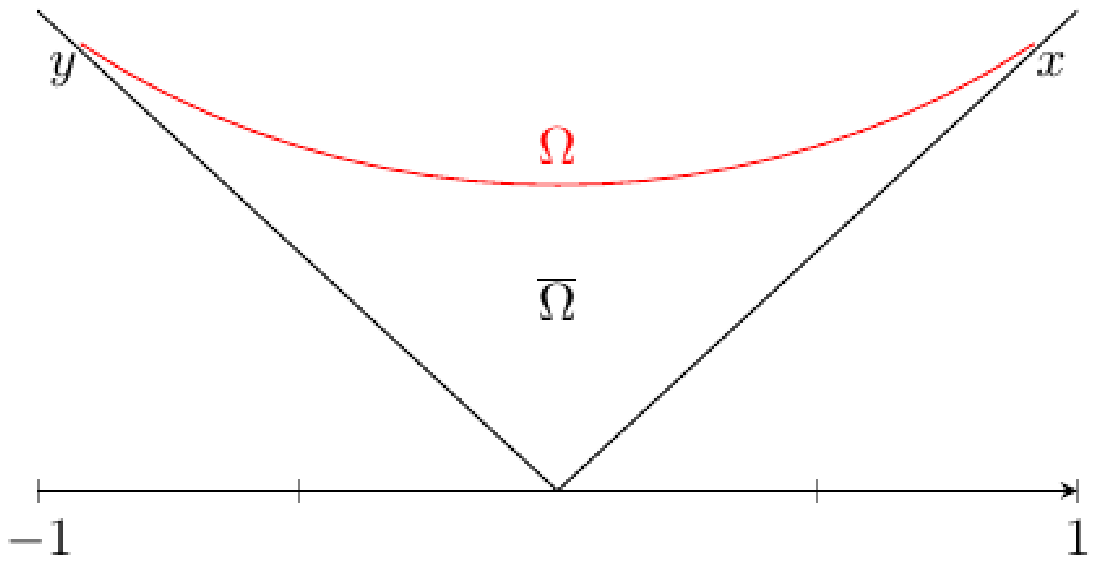}
\end{center}
\caption{}
\label{fig3}
\end{figure}

Recall the theorem on the limit shape of typical (with respect to the Planche\-rel measure) Young diagrams (\cite{VK}, \cite{LS}). It says that, under the normalization~$n^{-1/2}$, the sequence of finite diagrams converges (in the Hausdorff metric with respect to the Euclidean norm)
to a curvilinear triangle~$\bar\Omega$ (see Fig.~\ref{fig3}). The curvilinear side of this triangle is the graph of the function
$$
\Omega(s)=\frac{2}{\pi} (s\cdot \arcsin s+\sqrt{1-s^2}),
$$
where $|s|\leq 1$.

We want to find the limit shape of  the $P$- and $Q$-tableaux after an appropriate normalization.

Consider the sequences of $P$- and $Q$-tableaux corresponding to a realization $\{x_1,x_2,\dots \}\in I$. We introduce sequences of step functions $\{\phi_n\}$ and $\{\psi_n\}$ that are defined on diagrams regarded as subsets of the lattice~$\mathbb{Z}_+^2$ and take constant values inside every open cell. The functions~$\phi_n$ replace the $P$-tableaux and take values in the interval~$[0,1]$ according to the numbers occupying their cells; the functions~$\psi_n$ take positive integer values according to the numbers occupying the cells of the $Q$-tableaux.

Now we introduce the following scaling: for every~$n$, shrink the domain of definition, i.e., the diagram, by the factor~${\sqrt n}$, leaving the values of the first function unchanged and dividing the values of the second function  by~$n$. For the new step functions, defined on some subset of the plane~$\mathbb{R}^2_+$, we keep the previous notation: $\phi_n$~and~$\psi_n$. These two sequences are quite different in nature: before the normalization, the first of them had a zero weak limit, while the second one stabilized.

It follows from the theorem on the limit shape of Young diagrams with respect to the Plancherel measure that the body of the diagram,  under an appropriate normalization, turns into the curvilinear triangle~$\bar \Omega$, which can be thought of as a~part of the plane~$\mathbb{R}^2_+$. We introduce on $\bar \Omega$ polar coordinates $(r,\theta)$, where $\theta \in (0,\pi/2)$, $r \in (0,r_{\theta})$, and $r_{\theta}$ is the length of the radius vector on the ray corresponding to the angle~$\theta$ at a point of~$\Omega$.
Consider the following function $A(r,\theta)$ on the curvilinear triangle $\bar \Omega$ defined in polar coordinates:
$$
A(r,\theta)=r^2\!/r_{\theta}^2,\qquad A(r_{\theta},\theta)=1.
$$

One may say that the graph of this function is obtained by lifting the curve~$\Omega$ to height~$1$ with subsequent parabolic interpolation between~$0$ and the intersection point of the ray corresponding to the angle~$\theta$ with~$\Omega$ in each of the vertical planes.

\begin{figure}[h]
\begin{center}
\includegraphics[width=0.85\textwidth]{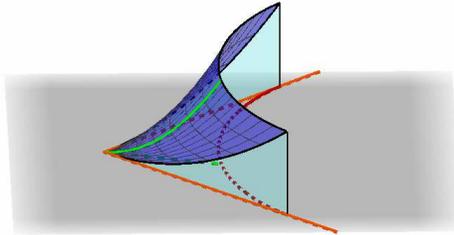}
\end{center}
\caption{The arch, i.e., the 3D limit shape of both Young tableaux}
\label{fig4}
\end{figure}

This function $A(\bcdot,\bcdot)$ will be called the \textit{arch}, according to the shape of its graph. It is depicted in Fig.~\ref{fig4} made by P.~B.~Zatitskii.

\begin{theorem}
The limits, in the measure~$m$, of the sequences of random step  functions $\{\phi_n\}_n$ and $\{\psi_n\}_n$ as $n\to \infty$ coincide and are equal to the function~$A$:
$$
\lim_n \phi_n(\bcdot,\bcdot)=A(\bcdot,\bcdot), \qquad \lim_n \psi_n(\bcdot,\bcdot)=A(\bcdot,\bcdot).
$$
\end{theorem}

\begin{pf}{Proof}
First, we refine the statement: the piecewise constant functions $\phi_n, \psi_n$ are defined as functions on the plane, but these functions are random, since they depend on the realization of the Bernoulli scheme. By convergence we mean that, in the measure~$m$ (or even almost surely), they converge to the function~$A$ in the uniform metric.

We emphasize once again that these two sequences of functions are very different in nature, however, their limits coincide. Accordingly, the proofs are based on quite different considerations. As to the first equality (for the functions~$\phi_n$, i.e., for the $Q$-tableaux), it is an easy corollary of the main theorem on the limit shape of Young diagrams with respect to the Plancherel measure (\cite{VK}, \cite{LS}), which was proved in~\cite{Gr}; it suffices to observe that the limit shape is homogeneous of degree~$2$ with respect to homotheties of the lattice.

To the author's knowledge, the second equality has not been considered earlier, as well as the very dynamics of $P$-tableaux in this  context.

As mentioned above, it is convenient to work with the functions~$\psi_n$ in polar coordinates and to restrict each of them to the rays
$L_{\theta}=(0,r_{\theta})$. This means that we choose a direction and consider the cells of the diagram intersected by the corresponding ray; look at the numbers occupying these cells.

We begin with corner cells. On every ray containing a corner cell, for almost all realizations, the numbers occupying this cell monotonically increase to~$1$; in other words, the limits of~$\psi_n$ at the corner cells are equal to~$1$ and thus coincide with the corresponding values of the function~$A$. Indeed, the monotone convergence is obvious, and the assumption that the limit is strictly less than~$1$ would contradict the density of the values of the realization.

Using the homothety with coefficient~$a$, where $|a|<1$,  centered at the zero vertex of the curvilinear triangle $\bar \Omega$,
we can apply the same reasoning to the image
$\Omega_a$ of the curve~$\Omega$ under this homothety, and thus extend the arguments to the whole triangle~$\bar \Omega$.
\qed\end{pf}

It would be extremely interesting to study the behavior of normalized fluctuations of values in cells: consider the random function
$$
(1-\psi_n(r_{\theta}))/\sqrt n
$$
 on the curve~$\Omega$. What is the limiting distribution of this function?

On the other hand, we can consider fluctuations on rays starting at the zero vertex,
i.e., the similar normalized differences of values of the functions $A$~and~$\psi_n$  on rays. These fluctuations
form an interesting random field on the triangle~$\bar \Omega$.

These considerations should correlate with theorems on fluctuations, for example, of the length of the first row of Young diagrams in the spirit of~\cite{BDJ},~\cite{BOO}.

Note also that all results of this section on the asymptotic behavior of $P$- and $Q$-tableaux remain valid  for the limiting behavior of these tableaux for permutations, i.e., are applicable to the asymptotics of the symmetric groups~$S_n$; the Bernoulli scheme in these considerations is irrelevant.

\subsection{Conclusion}
Above, we considered the structure of a Bernoulli scheme, i.e., a sequence of independent random variables with values in a linearly ordered set, that is, the interrelations between two orders: the given linear order and the temporal order.

It is of interest how the results of this and other related papers change if we

(a) change the distribution of independent variables, but preserve the range of values; in this case, the RSK correspondence still makes sense, and the question is about the behavior of $P$-tableaux;

(b) drop the independence; it would be interesting here to obtain results on a~combinatorial encoding similar to the one considered in this paper, and study the behavior of the corresponding $Q$-tableaux;

(c) drop the linear order on the space of values and consider other partial orders or even other structures on this space; the simplest example is when the values of random variables lie in a finite-dimensional cube; in this case, it is not even known what is the analog of the RSK correspondence.

On the other hand, for the classical Bernoulli scheme considered in this paper, the question about other invariant combinatorial encodings that are decodable (i.e., have the distinguishability property) is open.

\end{document}